\documentclass{article}
   \usepackage{amsmath,amssymb,amsthm,amscd}

   \newtheorem{Thm}{Theorem}[section]

   \newtheorem{Def}[Thm]{Definition}
   \newtheorem{Rem}[Thm]{Remark}
   \newtheorem{Exam}[Thm]{Example}
   
   \makeatletter
   \@addtoreset{equation}{section}
   \makeatother
   
\newcommand{\R}{\mathbb{R}}
\newcommand{\N}{\mathbb{N}}
\newcommand{\Z}{\mathbb{Z}}
\newcommand{\C}{\mathbb{C}}

\newcommand{\dbar}{{d\hspace{-0,05cm}\bar{}\hspace{0,05cm}}}
\newcommand{\Op}{\textup{Op}}
\newcommand{\op}{\textup{op}}

\newcommand{\g}{\textbf{g}}

\newcommand{\ci}{C^\infty}

   \author{B.-W. Schulze}
\title{The Iterative Structure of Corner Operators}

 \date{}
\begin{document}
  \maketitle

\begin{abstract}
\noindent
We give a brief survey on some new developments on elliptic operators on manifolds with polyhedral singularities.
The material essentially corresponds to a talk given by the author during the Conference ``Elliptic
and Hyperbolic Equations on Singular Spaces", October 27 - 31, 2008, at the MSRI, University of Berkeley. 

\end{abstract}

\noindent {\bf Mathematics Subject Classification: } \\
\noindent {\bf Primary: 35S35}\\
\noindent {\bf Secondary: 35J70} 

\tableofcontents
\section*{Introduction}
\addcontentsline{toc}{section}{Introduction}
\markboth{INTRODUCTION}{INTRODUCTION}
Manifolds with higher corners or edges of order  $k \in \N$ are (in our notation) special stratified spaces, where $k=0$ corresponds to smoothness, $k = 1$ to conical or edge singularities. Manifolds with singularities of order $k$ form a category $ \mathfrak{M}_k.$ Ellipticity of operators will be expressed by a principal symbolic hierarchy
$$\sigma = (\sigma_j)_{0 \leq j \leq k}$$
with $\sigma_0$ being the standard homogeneous principal symbol on the main stratum $\textup{int}\,M$, while the components $\sigma_j,\,j>0,$ live on the other strata and are operator-valued.
\begin{Exam}
 \begin{description}
  \item \textup{(i)}The half-axis $\overline{\R}_+$ can be regarded as a manifold with conical singularity $0$. 
  \item \textup{(ii)}The half-space $\overline{\R}_+  \times \Omega$ for an open set $\Omega \subseteq  \R^q$ is an example of a manifold with edge $\Omega $
and model cone $\overline{\R}_+$.
  \item \textup{(iii)}Let $X$ be a closed compact $C^{\infty}$ manifold, then the quotient space $X^{\Delta} := (\overline{\R}_+ \times X)/(\{0\} \times X)$ which is an infinite cone with vertex $v$, represented by $\{0\} \times X$, is a manifold with conical singularity $v$ and base $X$.
  \item \textup{(iv)}The wedge $X^{\Delta} \times \Omega$ with $X$ and $\Omega$ as before is a manifold with edge $\Omega$ and model cone $X^{\Delta}$. 
 \end{description}
\end{Exam}
\begin{Rem}
Consider a Riemannian metric of the form $dr^2 + r^2g_X$ on the open stretched cone $X^\wedge := \R_+\times X$ where $g_X$ is a Riemannian metric on $X$. Then for the associated Laplace-Beltrami operator we obtain \textup{(}for $m=2$\textup{)}
\begin{equation}
\label{i.1}
A = r^{-m}\sum_{j=0}^{m}a_j(r)\left( -r\frac{\partial }{\partial r}\right) ^j,
\end {equation}
with coefficients $a_j\in \ci (\overline{\R}_+, \textup{Diff}^{m -j}(X))$.
More generally, if we consider  a Riemannian metric $dr^2 + r^2g_X + dy^2$ on the open stretched wedge $X^\wedge \times \Omega$, then the associated Laplace-Beltrami operator has the form \textup{(}for $m=2$\textup{)}
\begin{equation}
\label{i.2}
A=r^{-m}\sum_{j+|\alpha |\leq m }a_{j\alpha }(r,y)\left( -r\frac{\partial }{\partial r}\right)^j(rD_y)^\alpha .
\end {equation}
with coefficients $a_{j\alpha }\in \ci(\overline{\R}_+\times \Omega , \textup{Diff}^{m -(j+|\alpha|)}(X))$.
\end{Rem}
The principal symbolic hierarchies are as follows.
In the conical case we have $\sigma_0(A) \in \ci(T^*X \setminus \{0\})$, the usual homogeneous principal symbol of $A$ (degenerate at $r=0$), and
\begin{equation}
\label{i.3}
\sigma _1(A)(z) = \sum_{j=0}^{m }a_j(0)z^j : H^s(X) \rightarrow H^{s-m}(X),
\end {equation}
the principal conormal symbol. In the edge case we have $\sigma_0(A) \in \ci(T^*(X \times \Omega)\setminus \{0\})$, the usual homogeneous principal symbol of $A$ (here edge-degenerate at $r=0$), and 
\begin{equation}
\label{i.4}
\sigma _1(A)(y,\eta)=r^{-m}\sum_{j+|\alpha |\leq m }a_{j\alpha }(0,y)\left( -r\frac{\partial }{\partial r}\right)^j(r \eta)^\alpha : \mathcal{K}^{s,\gamma}(X^\wedge)
\rightarrow \mathcal{K}^{s-m,\gamma-m}(X^\wedge),
\end {equation}
the homogeneous principal edge symbol (the meaning of $\mathcal{K}^{s,\gamma}(X^\wedge)$ will be explained below; the notation $\mathcal{K}$ comes from ``Kegel"). Note that there is a similarity between edge symbols and boundary symbols of differential operators on a manifold with
smooth boundary. Consider, for instance, a differential operator $A = \sum_{k+ \beta \leq m}b_{k \beta}(r,y)D_r^kD_y^{\alpha}$ on the half-space, $\overline{\R}_+  \times \Omega$ with coefficients $b_{k \beta}\in \ci(\overline{\R}_+  \times \Omega)$. Apart from 
\begin{equation}
\sigma _0(A)(r,y,\varrho,\eta)=\sum_{k+|\beta |= m }b_{k\beta }(r,y) \varrho^k \eta^{\beta}
\end {equation}
we have the principal boundary symbol
\begin{equation}
\sigma _1(A)(y,\eta)=\sum_{k+|\beta|= m}b_{k\beta }(0,y)D_r^k \eta^\beta : H^s(\R_+) \rightarrow H^{s-m}(\R_+),\quad s \in \R.
\end {equation}
Boundary symbols are homogeneous in the sense
$$\sigma_1(A)(y,\lambda \eta)= \lambda^\mu \kappa_{\lambda}\sigma_1(A)(y, \eta) \kappa_{\lambda}^{-1} \quad \textup{for all} \quad \lambda \in \R_+.$$
Here $\kappa_{\lambda} : H^s(\R_+) \rightarrow H^{s-m}(\R_+) $ is a strongly continuous group of isomorphisms, defined by $(\kappa_{\lambda}u)(r):= \lambda^{1/2}u(\lambda r), \lambda \in \R_+$. A similar relation holds for edge symbols, based on
$$\kappa_{\lambda} : \mathcal{K}^{s,\gamma}(X^ \wedge) \rightarrow \mathcal{K}^{s- m, \gamma - m}(X^ \wedge), s, \gamma \in \R, $$ where $(\kappa_{\lambda}u)(r,x):=
\lambda^{(n+1)/2}u(\lambda r,x),\, \lambda \in \R_+.$\\
In this presentation we give an idea on how to formulate algebras of (pseudo-differential) operators on $\textup{int}\,M$ that contain the (for the nature of singularities typical) differential operators, together with the parametrices of elliptic elements. More details may be found in \cite{Schu25}, \cite{Schu27}, and in a new monograph in preparation \cite{Schu58}, see also the references below.

% Section 1
\section{The category $\mathfrak{M}_k$}
Stratified spaces of different kind occur in numerous fields of mathematics and also in the applied sciences. Here we single out specific categories of such spaces where certain elements of the analysis of PDE can be formulated in an iterative manner. General references on stratified spaces are  Fulton and MacPherson \cite{Fult1}, or  Weinberger \cite{Wein1}.
\begin{Def}
A topological space $M$ \textup{(}under some natural conditions on the topology in general\textup{)} is said to be a manifold with  singularities of order $k \in \N, k \geq 1,$ if
 \begin{description}
  \item \textup{(i)} $M$ contains a subspace $Y \in \mathfrak{M}_0$ such that $M \setminus Y \in \mathfrak{M}_{k-1}$;

  \item \textup{(ii)}  $Y$ has a neighbourhood $V \subseteq M$ which is a \textup{(}locally trivial\textup{)} cone bundle over $Y$ with fibre
$X^\Delta$ for some $X \in \mathfrak{M}_{k-1}$.
 \end{description}
\end{Def}
Transition maps $X^\Delta  \rightarrow X^\Delta$ are induced by restrictions of $\mathfrak{M}_{k-1}$ -isomorphisms $ \R \times X \rightarrow \R \times X$ to $\overline{\R}_+\times X$. This gives rise to corresponding transition maps $X^\Delta  \times \Omega \rightarrow X^\Delta \times \tilde{\Omega}$ for the respective  $X^\Delta  $ -bundles over $Y$.

\begin{Rem}
$\mathfrak{M}_k$ is a category with a natural notion of morphisms and isomorphisms.
\end{Rem}
\begin{Rem}
$Y=: Y^k$ is called the minimal stratum of $M$. The space $M \setminus Y^k \in \mathfrak{M}_{k-1}$ contains a space $Y^{k-1} \in \frak{M}_0$, such that $(M \setminus Y^k) \setminus Y^{k-1} \in \frak{M}_{k-2}$, etc. This yields a representation
$$M= Y^0\cup Y^1\cup Y^2\cup ...\cup Y^k$$ as a disjoint union of strata $Y^j \in \mathfrak{M}_0$. We set $\textup{int} \, M := Y^0$, called the maximal stratum of M, and $\textup{dim} \, M := \textup{dim}(\textup{int} \, M)$. Moreover, $M$ is locally near $Y^j$ modelled on an $X_{j-1}^\Delta $-bundle over $Y^j$, for $X_{j-1}^\Delta \in \mathfrak{M}_{j-1}$ with trivialisations $X_{j-1}^\Delta \times \Omega^j, \Omega^j  \subseteq \R^{\textup{dim}Y^j}$ open.
\end{Rem}
\begin{Rem}
 The same topological space $M$ can be stratified in different ways. For instance, we have $M=\R^n  \in \mathfrak{M}_0$ but also $M \in \mathfrak{M}_1$ when we set
$Y^1= \{0\}, Y^0= \R^n \setminus \{0\}.$
\end{Rem}

\begin{Rem}
 $M\in \frak{M}_{k}$, $L\in \frak{M}_{l}$ implies $M \times L\in \frak{M}_{k+l}$.
\end{Rem}
There are many other interesting properties of the categories $\mathfrak{M}_k$ that we do not discuss in detail here. It would be desirable to develop the connection of our analysis on singular spaces with the work from topological side. For instance, D. Trotman informed me in Berkeley on his works with coauthors, cf. \cite{Bekk1} jointly with Bekka, and \cite{King1} with King.

\section{Corner-degenerate operators}
Let $\textup{Diff}_\textup{deg}^m(M)$ 
for $M \in \mathfrak{M}_k,  k \geq 1, $ denote the set of all
 $A \in \textup{Diff}_\textup{deg}^m(M \setminus Y)$,
($\textup{Diff}_{\textup{deg}}^m(\cdot) = \textup{Diff}^m(\cdot)$ in the smooth case) such that close to $Y$ in the local variables
$$(r,x)\in \R_+\times X \quad \textup{for}\quad \textup{dim} \, Y =0,\quad\textup{and}\quad
(r,x,y)\in \R_+\times X \times \Omega \quad \textup{for}\quad \textup{dim} \, Y > 0,$$
we have
$$A = r^{-m}\sum_{j=0}^{m}a_j(r)\left( -r\frac{\partial }{\partial r}\right) ^j \quad \textup{for}\quad \textup{dim} \, Y =0,$$
with coefficients $a_j\in \ci (\overline{\R}_+, \textup{Diff}_\textup{deg}^{m -j}(X))$, and
$$A=r^{-m}\sum_{j+|\alpha |\leq m }a_{j\alpha }(r,y)\left( -r \frac{\partial }{\partial r}\right)^j(rD_y)^\alpha \quad \textup{for}\quad \textup{dim} \, Y > 0, $$
with coefficients $a_{j\alpha }\in \ci(\overline{\R}_+\times \Omega , \textup{Diff}_\textup{deg}^{m -(j+|\alpha|)}(X))$, respectively.
The principal symbolic hierarchies are iteratively defined by
\begin{equation} \label{2.1}
\sigma (A) := (\sigma (A|_{M\setminus Y}), \sigma _k(A))
\end{equation}
where $\sigma (A|_{M\setminus Y})$ is known by the steps before, while
\begin{equation} \label{2.2}
\sigma _k(A)(z) = \sum_{j=0}^m a_j(0)z^j \quad \textup{for} \quad \textup{dim}\,Y =0, z \in \C, 
\end{equation}
and
\begin{equation} \label{2.3}
\sigma _k(A)(y, \eta)=r^{-m}\sum_{j+|\alpha |\leq m }a_{j\alpha }(0,y)\left( -r\frac{\partial }{\partial r}\right)^j(r \eta)^\alpha \quad \textup{for}
\quad \textup{dim} \, Y > 0, (y,\eta) \in T^*\Omega \setminus 0. 
\end{equation}
$\sigma _k(A)(z)$ takes values in $\textup{Diff}_\textup{deg}^m (X)$ for $\textup{dim}\,Y =0,$ and $\sigma _k(A)(y, \eta)$ takes values in  $\textup{Diff}_\textup{deg}^m (X^\wedge)$ for $\textup{dim} \, Y > 0.$
\begin{Rem}
If we dissolve the information in $\eqref{2.1}$ with respect to the other strata of $M$ we obtain $k+1$ components of $\sigma (A)$, namely,
$$\sigma (A) = (\sigma_0 (A), \sigma_1 (A), \ldots, \sigma _k(A)),$$
with $\sigma_0 (A) \in \ci(T^*(\textup{int} \, M) \setminus \{0\})$ being the standard homogeneous principal symbol of $A$ on the main stratum of $M$, while the other components $\sigma _j(A)$ are operator-valued and associated with $Y^j$, where $\sigma_j(A)$ is of analogous form as $\eqref{2.2}$ for $\textup{dim} \, Y =0,$ and $\eqref{2.3}$ for $\textup{dim} \, Y > 0$ \textup{(}in the latter case parametrised by $T^*Y^j \setminus 0$\textup{)}. The symbol $\sigma _j(A)$ for $\textup{dim} \, Y =0$ acts between
weighted spaces $H^{s,\gamma(j-1)}(X_{j-1})$ where M is locally near $Y^j$modelled on $X_{j-1}^\Delta$. Moreover, $\sigma _j(A)$ for $\textup{dim} \, Y > 0$ is analogous to $\eqref{2.3}$, now parametrised by points in $T^*Y^j \setminus 0$ and acting between spaces of the form $\mathcal{K}^{s,\gamma(j)}(X_{(j-1)}^ \wedge)$ where $M$ is modelled on $X_{j-1}^ \wedge \times \Omega^j$ locally near $Y^j,$ and $\gamma(j-1)$ is a $(j-1)$-tuple of weights.
\end{Rem}

\section{Problems and results}

The analysis of operators of the spaces $\textup{Diff}_\textup{deg}^m(M)$ on stratified spaces $M \in \mathfrak{M}_k,  k \in \N, $ gives rise to a number of \textbf{natural problems}
that are solved by works of several authors in this field or are open and still represent \textbf{challenges} for the future development. Let us give a list of such problems, and then some key words concerning results and references: \\ \\
$\bullet$ What is ellipticity of $A \in \textup{Diff}_\textup{deg}^m(M)$ in connection with the principal symbolic hierarchies explained in the preceding section? \\ \\
$\bullet$ Construct a pseudo-differential calculus containing $\textup{Diff}_\textup{deg}^m(M)$ together with the parametrices of elliptic elements. \\ \\
$\bullet$ Establish the Fredholm property and study the index of elliptic operators in weighted distribution spaces when $M$ is compact. \\ \\
$\bullet$ Understand ellipticity, parametrices, and Fredholm property in suitable weighted spaces when $M$ has conical exits to infinity. \\ \\
$\bullet$ Study parameter-dependent theories on $M$. \\ \\
$\bullet$ Characterise asymptotics of solutions to elliptic equations (discrete, continuous, variable/branching, or iterated), in simplest cases of the form
\begin{equation} \label{3.1}
u(r, \cdot) \sim \sum_j \sum_{k=0}^{m_j} c_{jk}(\cdot) r^{-p_j} \textup{log}^kr \quad \textup{as} \quad r \rightarrow 0,
\end{equation}
with $p_j \in \C, \textup{Re}\, p_j \rightarrow -\infty \quad \textup{as} \,r \rightarrow \infty$ (if the expansion is infinite).\\ \\
$\bullet$ Compute the points $p_j$ for interesting examples; those points appear as ``non-linear eigenvalues" of conormal symbols.\\ \\ 
$\bullet$ Study various quantisations of corner-degenerate symbols, in particular, in terms of holomorphic/meromorphic operator functions with values in algebras of lower singularity order.\\ \\
$\bullet$ Understand the hierarchy of topological obstructions appearing in the construction of elliptic operators with prescribed elliptic symbols $\sigma_0$.\\ \\
$\bullet$ Study index theories, homotopy classifications, K\"unneth formulas, etc., for the higher corner operator algebras.\\ \\
\textbf{Models from diverse applications with singular geometry: }\\ \\
$\bullet$ Mixed problems, operators with/without the transmission property at the boundary (Zaremba problem, etc.). \\ \\
$\bullet$ Boundary value problems in polyhedral domains with the induced metric from an ambient space, occurring in elasticity, mechanics (beams, shells, plates, \ldots). \\\\
$\bullet$ Crack problems with crack boundaries that are smooth or have singularities. \\ \\
$\bullet$ Operators with singular potentials
$$\Delta + V, \quad V=\sum_{i,j=0}^{N} c_{ij}|x^{(i)}-x^{(j)}|^{-1},$$ for the Laplacian $\Delta$ in $\R^{3N}$, and $x^{(j)}= (x_1^{(j)},x_2^{(j)},x_3^{(j)})$ indicating
the position of three-dimensional particles; the question is to describe the behaviour of solutions to $(\Delta + V)u =f,$ close to the singularities of the potential $V,$ say, for smooth $f.$ \\ \\
\textbf{The operator algebras that we discuss here contain many special cases and substructures:} \\ \\
$\bullet$ Singular integral operators with piecewise smooth coefficients (cf. Gohberg and Krupnik \cite{Gohb2}). \\ \\
$\bullet$ Mellin operators on the half-axis (cf. Eskin \cite{Eski2}). \\ \\
$\bullet$ Operators on manifolds with conical exits to infinity (cf. Shubin \cite{Shub1}, Parenti \cite{Pare1}, Cordes \cite{Cord1}) \\ \\ 
$\bullet$ Parameter-dependent operators (cf. Agranovich and Vishik \cite{Agra1}).\\ \\ 
$\bullet$ Boundary value problems without/with the transmission property at the boundary (cf. Vishik and Eskin \cite{Vivs2}, Eskin \cite{Eski2}, Boutet de Monvel \cite{Bout1}). \\ \\
$\bullet$ Totally characteristic operators (cf. Melrose \cite{Melr4}, Melrose and Mendoza \cite{Melr1}). \\ \\
$\bullet$ Edge-degenerate, and corner-degenerate operators (cf. Rempel and Schulze \cite{Remp3}, Schulze \cite{Schu32}, Mazzeo \cite{Mazz1}, Schulze \cite{Schu29}). \\ \\
$\bullet$ Boundary value problems in the frame of the edge calculus (cf. Rempel and Schulze \cite{Remp1}, Schulze \cite{Schu31}, Schulze and Seiler \cite{Schu41}). \\ \\
Other authors of the research group of the University of Potsdam and guests during the past years contributed improvements and important new aspects, cf. Witt \cite{Witt3}, Gil \cite{Gil4}, Gil and Mendoza \cite{Gil5}, Seiler \cite{Seil1}, \cite{Seil2}, Krainer \cite{Krai9}, \cite{Krai5}, or the author's joint papers with
Coriasco \cite{Cori2}, Dines and Liu \cite{Dine3}, Flad and Schneider \cite{Flad1}, Wei \cite{Schu56}; other references are given below. \\ \\There also appeared (or are in preparation) some monographs on these topics, in particular, \cite{Schu2}, \cite{Schu31}, \cite{Schu20}, or, jointly with Egorov \cite{Egor1}, Kapanadze \cite{Kapa10}, Harutyunyan \cite{Haru13}, Volpato \cite{Schu60}. More details on the higher corner calculus will also be given in the author's new monograph \cite{Schu58}.\\ \\

\section{Some typical tools of the higher corner calculus}
Let us first consider $k=0$ which is the smooth case. On a $C^ \infty$ manifold $M$ we have $L^m_{\textup{cl}}(M)$, the space of classical pseudo-differential operators of order $m \in \R$ (``classical" is here not essential, but for $k>0$ we employ this assumption). For an $A \in L^m_{\textup{cl}}(M)$ we have the standard homogeneous principal symbol $\sigma_0(A) \in \ci(T^*M \setminus \{0\}).$ Clearly everything works for vector bundles as well. \\
Let $H_{\textup{comp}}^s(M), H_{\textup{loc}}^s(M), s \in \R,$ denote the standard Sobolev spaces, and write $H^s(M)$ when $M$ is compact or an Euclidean
 space. The spaces $H^s(\R^n \times \R^q)$ admit anisotropic reformulations that are useful for the singular cases, namely, $$H^s(\R^n\times \R^q) = \mathcal{W}^s(\R^q,H^s(\R^n)).$$ Here $\mathcal{W}^s(\R^q,H)$ for some Hilbert space $H$ which is endowed with a strongly continuous group $\{\kappa_{\lambda}\}_{\lambda \in \R_+}$ of isomorphisms $\kappa_{\lambda}: H \rightarrow H$ means the completion of $\mathcal{S}(\R^q, H)$ with respect to the norm $\{\int\langle \eta\rangle ^{2s} \|\kappa_{\langle \eta \rangle}^{-1} \hat{u}(\eta)\|_H^2 d\eta\}^{1/2}, 
\langle\eta\rangle= (1+|\eta|^2)^{1/2}.$ In the case $H=H^s(\R^n)$ we set $(\kappa_{\lambda} u)(x)= {\lambda}^{n/2}u(\lambda x)$.\\

Let $L_{\textup{cl}}^m(M,\R^l)$ denote the space of parameter-dependent pseudo-differential operators with parameter $\lambda \in \R^l, l \in \N,$ on an open $C^\infty$ manifold $M$, with local amplitude functions $a(x,\xi,\lambda)$ that are classical symbols in $(\xi,\lambda) \in \R^{n+l}, n= \textup{dim}\, M$, and $L^{- \infty}(M,\R^l)= \mathcal{S}(\R^l,L^{-\infty}(M)),$ where $L^{-\infty}(M)$ is identified with $\ci(M \times\ M)$ via a fixed Riemannian metric. We employ the fact that for compact $M$ there exist parameter-dependent elliptic order reducing isomorphisms
$$R^m(\lambda):H^s(M) \rightarrow H^{s-m}(M)$$ for every $m$ and $s$.
Let us now give an idea on how the respective parameter-dependent operator spaces $\mathfrak{A}^m(M,\g;\R^l),m \in \R,$ are constructed in the case $M \in \mathfrak{M}_1$ which corresponds to conical or edge singularities (and also contains the case of a manifold with smooth boundary). Here $\lambda \in \R^l$ is the parameter, and $\g= (\gamma, \gamma - m, \Theta)$ are weight data for a weight $\gamma \in \R $ and a weight interval $\Theta = (\theta, 0]$ for a $ -\infty \leq \theta < 0$, where we control asymptotics.\\
Let us forget about $\R^l$ for a while, and define weighted spaces, first for conical singularities. By $\mathcal{H}^{s,\gamma}(X^\wedge)$ on the open stretched cone $X^\wedge = \R_+\times X$ we denote the completion of $\ci_0 (X^\wedge)$ with respect to the norm
$${\Big {\{}}(2\pi i)^{-1} \int_{\Gamma_{{\frac{n+1}{2}}-\gamma }} \|R^s(\textup{Im}z){Mu(z)\|_{L^2(X)}^2 dz\Big {\}}}^{1/2}$$
where $R^s(\lambda) \in L_{\textup{cl}}^s(X;\R)$ is a parameter-dependent elliptic family, $n=\textup{dim}\,X,\, \Gamma_\beta = \{z \in \C: \textup{Im}\,z = \beta\},$ and $Mu(z)= \int_0^\infty r^{z-1}u(r)dr $ is the Mellin transform on $\R_+.$ The space $\mathcal{H}^{s,\gamma}(X^\wedge)$ takes part in the definition of the space $\mathcal{K}^{s,\gamma}(X^\wedge),$ namely, close to $r=0.$ Another ingredient close to $r=\infty$ is the space $H^s_{\textup{cone}}(X^\wedge).$ Let us first define a version on $\R \times X $ rather than $X^\wedge.$ The space $H^s_{\textup{cone}}(\R \times X)$ is defined to be the completion of $\ci (\R \times X)$ with respect to the norm $$\Big \{ \int \|\langle r\rangle ^{-s}\textup{Op}_r(p)(\eta^1)u\|_{L^2(X)}^2 dr\Big \}^{1/2},$$ $\langle  r \rangle = (1+r^2)^{1/2}$, $p(r,\rho,\eta)= \tilde{p}(\langle r\rangle \rho,\langle r \rangle \eta), $ where $\tilde{p}(\tilde \rho,\tilde \eta)\in L_{\textup{cl}}^s(X; \R^{1+q}_{\tilde \rho,\tilde \eta})$ is a parameter-dependent elliptic family on $X,$ and $|\eta^1|$ is sufficiently large and fixed. Moreover, $$\textup{Op}_r(p)u(r) = \iint e^{i(r-r')\rho}p(r,\rho)u(r')dr'\dbar\rho,$$ 
$\dbar\rho := (2\pi)^{-1}d\rho.$ We set $H^s_{\textup{cone}}(X^\wedge):= H^s_{\textup{cone}}(\R \times X)|_{X^\wedge}.$ For any cut-off function $\omega(r)$ we define
$$\mathcal{K}^{s,\gamma}(X^\wedge)= \omega\mathcal{H}^{s,\gamma}(X^\wedge) + (1-\omega)H^s_{\textup{cone}}(X^\wedge),$$
and $\mathcal{K}^{s,\gamma;g}(X^\wedge) := \langle r \rangle ^{-g}\mathcal{K}^{s,\gamma}(X^\wedge)$, $s,\gamma, g \in \R.$ Observe that we have $$H^s_{\textup{comp}}(\R^n\times\R^q)\subseteq \mathcal{W}^s(\R^q,\mathcal{K}^{s,\gamma}(X^\wedge))\subseteq H^s_{\textup{loc}}(\R^n\times\R^q)$$ for all $s,g\in\R.$
The pseudo-differential background of cone and edge operator algebras are degenerate operators of the form $r^{-m}\textup{Op}_{r,y}(p),$ based on the Fourier transform, for   $p(r,y,\rho,\eta)= \tilde{p}(r,y,r\rho,r\eta),$ where $\tilde{p}(r,y,\tilde \rho,\tilde \eta)\in \ci (\overline{\R}_+ \times \Omega ,L_{\textup{cl}}^m(X; \R^{1+q}_{\tilde \rho,\tilde \eta})).$ It is useful to pass to Mellin quantisations, i.e., to operator-valued symbols referring to the Mellin transform. Let us explain here the edge case, i.e., $q>0$ (the conical case is simpler). To this end we define the space $$M_{\mathcal{O}}^m(X;\R^q_\eta),$$ consisting of all $h(z,\eta) \in \mathcal{A}(\C_z,L_{\textup{cl}}^m(X;\R^q_\eta))$ such that $h(\beta + i\rho,\eta)\in L_{\textup{cl}}^\mu(X;\R^{1+q}_{\rho,\eta})$ for every $\beta \in \R,$ uniformly in finite $\beta$-intervals. Here $\mathcal{A}(U,E)$ for an open set $U\subseteq\C$ and a Fr\'echet space $E$ is the space of all holomorphic functions in $U$ with values in $E,$ in the topology of uniform convergence on compact subsets.\\
An inversion process in the construction of parametrices of elliptic operators gives rise to symbols of the kind $$f(y,z)\in \ci(\Omega,M^{-\infty}_R(X)).$$ Here $R$ is an asymptotic type, in the most precise version $y$-wise discrete, otherwise a continuous asymptotic type.
Let us give an idea of the discrete case. Then $R$ is a sequence $\{(r_j,n_j)\}_{j\in\Z}\subset \C\times\N$ with $|\textup{Re}\,r_j|\rightarrow\infty$
as $|j|\rightarrow\infty.$ The space $M^{-\infty}_R(X)$ is defined to be the set of all $f\in\mathcal{A}(\C\setminus\pi_\C R,L^{-\infty}(X)), \pi_\C R := \{r_j\}_{j\in\Z},$ such that $f$ is meromorphic with poles at the points $r_j$ of mutliplicity $n_j+1.$ Moreover, $f(z)$ is stongly decreasing as $|\textup{Im}\,z|\rightarrow\infty,$ i.e., if $\chi(z)$ is any $\pi_\C R$-excision function ($=0$ close to $\pi_\C R,$ and $=1$ when $\textup{dist}\,(z,\pi_\C R)>c$ for some $c>0$) then $\chi(z)f(z)|_{\Gamma_\beta}\in \mathcal{S}(\Gamma_\beta,L^{-\infty}(X))$ for every $\beta\in\R,$ uniformly in compact $\beta$-intervals; here $$\Gamma_\beta:= \{z\in\C:\textup{Re}\,(z)=\beta\}.$$ More generally we also employ so-called continuous asymptotic types $R,$ representd by sets $V\subset\C$ such that $V \cap \{c\leq\textup{Re}\,(z)\leq c'\}$ is compact for every $c\leq c',$ cf. \cite{Schu2}, \cite{Schu20}; $V:=\pi_\C R.$
\begin{Thm} \label{4mel}
For every $p(r,y,\rho,\eta) = \tilde{p}(r,y,r\rho,r\eta),$\,$\tilde{p}(r,y,\tilde{\rho},\tilde{\eta})\in \ci(\overline{\R}_+ \times \Omega,L^m (X ;\R^{1+q}_{\tilde{\rho},\tilde{\eta}}))$ there exists an $h(r,y,z,\eta)=\tilde{h}(r,y,z,r\eta)$ of the form $\tilde{h}(r,y,z,\tilde{\eta})\in \ci(\overline{\R}_+ \times \Omega,M_{\mathcal{O}}^m(X;\R^q_{\tilde{\eta}}))$ such that 
\begin{equation} \label{4.3}
\Op_{r,y}(p)=\op_M^{\gamma - n/2} \Op_y(h)\quad\textup{mod}\,L^{-\infty}(X^\wedge \times \Omega)
\end{equation}
 for every $\gamma  \in \R$.
\end{Thm}
The operator-valued amplitude functions
\begin{equation} \label{4.3a}
a(y,\eta) = \omega(r)r^{-m}\op_{M_r}^{\gamma - n/2}(h)(y,\eta)\tilde{\omega}(r) +\textup{ smoothing Mellin plus Green symbols}(y,\eta) 
\end{equation}
with cut-off functions $\omega(r),\tilde{\omega}(r)$ furnish the symbols of the edge pseudo-differential calculus near $r=0.$ Those are symbols as follows.
Let $H$ and $\tilde{H}$ be Hilbert spaces with group actions $\{\kappa_\lambda\}_{\lambda\in \R_+}$ and $\{\tilde{\kappa}_\lambda\}_{\lambda\in \R_+},$ respectively. Then
$$S^m(\Omega \times \R^q;H,\tilde{H})$$ 
for $m\in\R$ and  $\Omega\subseteq \R^p$  open is defined to be the set of all $a(y,\eta)\in \ci(\Omega\times\R^q,\mathcal{L}(H,\tilde{H}))$ such that
$$\|\tilde{\kappa}_{\langle\eta\rangle}^{-1}\{D_y^\alpha  D_\eta^\beta a(y,\eta)\}\kappa_{\langle \eta \rangle} \|_{\mathcal{L}(H,\tilde{H})} \leq c\langle\eta\rangle^{m-|\beta|},$$ 
uniformly on compact subsets of $\Omega,$ for all $\eta\in\R^q$ and all multi-indices $\alpha,\beta.$ The subspace $S^m_{\textup{cl}}(\Omega \times \R^q;H,\tilde{H})$ of classical symbols is defined in terms of asymptotic expansions $\sum_{j=0}^{\infty}\chi(\eta)a_{(\mu-j)}(y,\eta),$ where $\chi(\eta)$ is an excision function, and $a_{(m-j)}(y,\eta)\in \ci(\Omega\times (\R^q\setminus\{0\}),\mathcal{L}(H,\tilde{H}))$ are of twisted homogeneity $\mu-j,$ i.e., $$a_{(m-j)}(y,\lambda\eta)= \lambda^{\mu-j}\tilde{\kappa}_\lambda a_{(m-j)}(y,\eta)\kappa_\lambda^{-1},\,\lambda \in \R_+.$$ Parallel to such operator-valued symbols we have vector-valued spaces $\mathcal{W}^s(\R^q,H)$ for a Hilbert space $H$ with group action $\{\kappa_\lambda\}_{\lambda\in \R_+},$ defined as the completion of $\mathcal{S}(\R^q,H)$ with respect to the norm $\big{\{}\int \langle\eta\rangle^{2s}\|\kappa_{\langle\eta\rangle}^{-1}\hat{u}(\eta)\|_H^{1/2}d\eta\big{\}}^{1/2}.$ There is also a straightforward generalisation to spaces of the kind $\mathcal{W}_{\textup{comp}}^s(\Omega,H)$ and $\mathcal{W}_{\textup{loc}}^s(\Omega,H),$ respectively, over an open set $\Omega \subseteq \R^q.$
\begin{Thm}
The above-mentioned operator functions $a(y,\eta)$ of the form $\eqref{4.3a}$ belong to 
$$S^m(\Omega\times\R^q;\mathcal{K}^{s,\gamma}(X^\wedge),\mathcal{K}^{s-m,\gamma -m}(X^\wedge))$$ 
based on $\{\kappa_{\lambda}\}_{\lambda \in \R_+},$ defined by $(\kappa_{\lambda}u)(r,x)= \lambda^{(n+1)/2}u(\lambda r,x)$ for $u\in \mathcal{K}^{s,\gamma}(X^\wedge)$ and induce continuous operators
\begin{equation} \label{4.4}
\Op_y(a): \mathcal{W}_{\textup{comp}}^s(\Omega,\mathcal{K}^{s,\gamma}(X^\wedge)) \rightarrow \mathcal{W}_{\textup{loc}}^{s-m}(\Omega,\mathcal{K}^{s-m,\gamma-m}(X^\wedge)) 
\end{equation}
 for all $s \in \R$.
\end{Thm}
\begin{Rem}
Observe that Theorem $\ref{4mel}$ and the continuity $\eqref{4.4}$ show that $$r^{-m}p(r,y,\rho,\eta)\rightarrow a(y,\eta)\rightarrow \textup{Op}_y(a)$$represents an operator convention \textup{(}quantisation\textup{)} for edge-degenerate symbols $r^{-m}p$. In the author's joint paper \textup{\cite{Gil2}} with Gil and Seiler it has been proved that the first non-smoothing term in $\eqref{4.3a}$ is equivalent \textup{(}\textup{mod} Green operators\textup{)} to another earlier quantisation of \textup{\cite{Schu32}}.
\end{Rem}
Let us now define smoothing Mellin plus Green symbols, already occurring in $\eqref{4.3a}.$
\begin{Def}
\begin{description}
  \item \textup{(i)} A smoothing Mellin symbol is an element $g_M$ in $$\bigcap_{s\in\R}S^m_{\textup{cl}}(\Omega\times\R^q;\mathcal{K}^{s,\gamma}(X^\wedge),\mathcal{K}^{\infty,\gamma -m}(X^\wedge))$$ which has an asymptotic expansion \textup{mod}\,
$\bigcap_{s\in\R}S^\infty(\Omega\times\R^q;\mathcal{K}^{s,\gamma}(X^\wedge),\mathcal{K}^{\infty,\infty}(X^\wedge))$ into symbols of the form $$r^{-m+j}\omega(r[\eta])\textup{op}_M^{\gamma_{j\alpha}-n/2}(f_{j\alpha})(y)\eta^\alpha \tilde{\omega}(r[\eta]);$$ here $\eta\rightarrow[\eta]$ is a strictly positive function in $\ci(\R^q)$ such that $[\eta]=|\eta|$ for $|\eta|>c$ for some $c>0,$ moreover, $$f_{j\alpha}(y,z)\in\ci(\Omega,M^{-\infty}_{R_{j\alpha}}(X)),$$ for some asymptotic types $R_{j\alpha},$ and $|\alpha|\leq j,
\,\pi_\C R_{j\alpha}\cap \Gamma_{{\frac{n+1}{2}-\gamma_{j\alpha}}}= \emptyset,\, \gamma_{j\alpha}\leq \gamma\leq \gamma_{j\alpha}+j.$

  \item \textup{(ii)} A Green symbol $g(y,\eta)$ is defined by $$g(y,\eta)\in\bigcap_{s,g\in\R}S^m_{\textup{cl}}(\Omega\times\R^q; \mathcal{K}^{s,\gamma;g}(X^\wedge),\mathcal{S}_P(X^\wedge)),$$ and $$g^*(y,\eta)\in\bigcap_{s,g\in\R}S^m_{\textup{cl}}(\Omega\times\R^q;\mathcal{K}^{s,-m+\gamma;g}(X^\wedge),\mathcal{S}_Q(X^\wedge)),$$ for continuous asymptotic types $P,$ and $Q,$ respectively. \textup{(}Concerning details on continuous asymptotic types, see, for instance, \textup{\cite{Schu2}}, or \textup{\cite{Schu20}}.\textup{)}
 \end{description}
\end{Def}
For a smoothing Mellin symbol $g_M(y,\eta)$ we set $$\sigma_1(g_M)(y,\eta):= r^{-m}\omega(r|\eta|)\textup{op}_M^{\gamma_{00}-n/2}(f_{00})(y)\tilde{\omega}(r|\eta|)$$ which is the homogeneous principal part of order $m$ of the respective classical operator-valued symbol. Analogously, if $g(y,\eta)$ is a Green symbol we set $$\sigma_1(g)(y,\eta)= g_{(m)}(y,\eta)$$ with $g_{(m)}$ being the homogeneous principal part of $g$ of order $m.$ 
\begin{Rem}
Edge symbols $a(y,\eta)$ are an analogue of boundary symbols from boundary value problems for operators with/without the transmission property at the boundary.
\end{Rem}
Let $M$ be a manifold with edge $Y$ of dimension $q>0$. Then the space of all edge pseudo-differential operators on $M$, referring to the weight data $\g=(\gamma,\gamma -m,\Theta)$ for a weight $\gamma\in\R$ and a weight interval $\Theta =(\vartheta,0], -\infty\leq\vartheta <0$ (which indicates an interval on the left of $\gamma,$ and $\gamma -m,$ respectively, where we control asymptotics) is defined to be the subset $$\mathfrak{A}^m(M,\g) \subset L^m_{\textup{cl}}(M\setminus Y)$$ of all operators $A$ that are locally near $Y$ of the form $A=\textup{Op}_y(a)\,\textup{mod}\,\mathfrak{A}^{-\infty}(M,\g)$ where $a(y,\eta)$ is an edge amplitude function $\eqref{4.3a}$ while $\mathfrak{A}^{-\infty}(M,\g)$ is defined by mapping properties to smooth functions with asymptotics. In order not to overload the explanations we omit some details on asymptotics; let us only note that control of asymptotics in terms of $\Theta,$ for instance, in the discrete case $\eqref{3.1}$ means that we observe exponents such that $\textup{Re}\,p_j$ belong to the interval $((n+1)/2-\gamma -\vartheta,(n+1)/2-\gamma)$ for functions in the preimage and to $((n+1)/2-\gamma -m -\vartheta,(n+1)/2-\gamma -m)$ in the image. For a first understanding it suffices to imagine $\Theta = (-\infty,0];$ then we may forget about $\Theta$ and write $\g=(\gamma,\gamma-m).$\\
The principal symbolic structure of operators $A\in \mathfrak{A}^m(M,\g),\, m\in \R,\, \g=(\gamma,\gamma -m, \Theta)$ is defined by $\sigma(A)= (\sigma_0(A),\sigma_1(A))$ with $\sigma_0(A)$ being the homogeneous principal symbol in the sense of $\mathfrak{A}^m(M,\g) \subset L^m_{\textup{cl}}(M\setminus Y),$ which is locally near $r=0$ of the form $\sigma_0(A)= r^{-m}\tilde{p}_{(m)}(r,x,y,r\rho,\xi,r\eta)$ where $\tilde{p}_{(m)}(r,x,y,\tilde{\rho},\xi,\tilde{\eta})$ is the homogeneous principal symbol of the above family in $\ci(\overline{\R}_+  \times \Omega,L^m_{\textup{cl}}(X;\R^{1+q}_{\tilde{\rho},\tilde{\eta}})).$ Moreover, for the case $q>0$ we have $$\sigma_1(A)(y,\eta)= r^{-m}\op_M^{\gamma - n/2}(h_0)(y,\eta)+\sigma_1(g_M+g)(y,\eta),\, (y,\eta)\in T^*\Omega\setminus 0,\,h_0(r,y,z,\eta)=\tilde{h}(0,y,z,r\eta),$$ which is a family of linear continuous operators
\begin{equation} \label{4.51}
\sigma_1(A)(y,\eta):\mathcal{K}^{s,\gamma}(X^\wedge)\rightarrow \mathcal{K}^{s-m,\gamma -m}(X^\wedge)
\end{equation} 
of homogeneity $$\sigma_1(A)(y,\lambda\eta)=\lambda^m\kappa_{\lambda}\sigma_1(A)(y,\eta)\kappa_{\lambda}^{-1},\,\lambda\in\R_+.$$

\begin{Thm}
Every $A\in\mathfrak{A}^m(M,\g),\,\g=(\gamma,\gamma -m, \Theta),\,M $ compact, induces continuous operators
\begin{equation} \label{4.5}
 A: \mathcal{W}^{s,\gamma}(M)\rightarrow \mathcal{W}^{s-m,\gamma-m}(M)
\end{equation}
 for all $s \in \R$. The operator $\eqref{4.5}$ is compact when $\sigma(A)=0.$
\end{Thm}
\begin{Def} \label{4.6def}
An operator  $A\in\mathfrak{A}^m(M,\g),\,\g=(\gamma,\gamma -m, \Theta),$ is said to be elliptic if
 \begin{description}
  \item \textup{(i)} A is elliptic as an operator in $L^m_{\textup{cl}}(\textup{int}M),$ and if in addition locally near $r=0$ the function $\tilde{p}_{(m)}(r,x,y,\tilde{\rho},\xi,\tilde{\eta})$ does not vanish for all $(\tilde{\rho},\xi,\tilde{\eta}) \neq 0 $ up to $r=0$;

  \item \textup{(ii)} the operators $\eqref{4.51}$   are bijective for all $(y,\eta)\in T^*\Omega\setminus 0.$
 \end{description}
\end{Def}
\begin{Rem}
The second condition of ellipticity concerning $\sigma_1$ is stronger than necessary. It suffices to impose the Fredholm property together with a $2\times2$ block matrix extension of $\sigma_1$ by extra trace and potential symbols to a family of isomorphisms. The extra symbols represent additional operators satisfying an analogue of the Shapiro-Lopatinskij condition, known from boundary value problems. Similarly as in the latter case this requires vanishing of a topological obstruction for $\sigma_0(A)$
\textup{(}concerning more details on that point, including the calculus when this topological obstruction does not vanish, see \textup{\cite{Schu42}}\textup{)}.

\end{Rem}
\begin{Thm}
An operator $A\in\mathfrak{A}^m(M,\g),\,\g=(\gamma,\gamma -m, \Theta),\,M $ compact, is elliptic with respect to $(\sigma_0,\sigma_1)$ if and only if $\eqref{4.5}$ is Fredholm for some fixed $s\in\R.$ In general the ellipticity of $A$ entails the existence of a parametrix in $\mathfrak{A}^{-m}(M,\g^{-1})$ belonging to $(\sigma_0^{-1},\sigma_1^{-1}).$
\end{Thm}
\begin{Rem}
Parameter-dependent operators of the class $\mathfrak{A}^m(M,\g;\R^l)$ are defined in an analogous manner as for $l=0$. There is then a notion of parameter-dependent ellipticity. If $\mathfrak{A}^m(M,\g;\R^l)$ is parameter-dependent elliptic, M compact, then
\begin{equation} \label{4.6}
A(\lambda):\mathcal{W}^{s,\gamma}(M) \rightarrow \mathcal{W}^{s-m,\gamma-m}(M)
\end{equation}
are isomorphisms for all $\lambda \in \R^l,\,|\lambda|$ sufficiently large, $s\in \R.$
\end{Rem}
Let $\mathfrak{A}^{m-1}(M,\g;\R^l) := \big\{A\in \mathfrak{A}^m(M,\g;\R^l): \sigma(A)=0\big\},$ and successively define $\mathfrak{A}^{m-j}(M,\g;\R^l)$ for every $j\in \N,\,\g=(\gamma,\gamma -m, \Theta).$
\begin{Thm}
For every $s', s''\in \R$ and $N\in \N$ there exists a $j\in \N$ such that for $A(\lambda) \in \mathfrak{A}^{m-j}(M,\g;\R^l)$ we have
\begin{equation} \label{4.7}
\|A(\lambda)\|_{\mathcal{L}(\mathcal{W}^{s',\gamma}(M),\mathcal{W}^{s'',\gamma-m}(M))} \leq c\langle\lambda\rangle^{-N}
\end{equation}
for all $\lambda\in\R $ and some $c>0.$
\end{Thm}

\section{Higher corner operators} 
We sketch a number of structures of the pseudo-differential operator calculus on a manifold with higher corners. For convenience we focus the consideration on the case $k=2.$ It will be fairly obvious that the concept is iterative and can be applied for higher corners as well, cf. \cite{Schu25}. This aspect is one of the main motivations of our approach.
Another motivation is, of course, to express parametrices of elliptic elements within the calculus which belongs to one of our results; special cases have been treated before, cf., \cite{Schu29}, \cite{Schu27}. Other contributions to the higher corner calculus are \cite{Schu25}, and the author's joint papers with Maniccia \cite{Mani2}, Krainer \cite{Krai3}, Calvo and Martin \cite{Calv2}, Calvo \cite{Calv3}, Harutyunyan \cite{Haru11}, \cite{Haru12}, \cite{Haru13}. 
Let $B\in\mathfrak{M}_1;$ then a starting point are corner-degenerate families
\begin{equation} \label{5.1}
p(t,z,\tau,\zeta)= \tilde{p}(t,z,t\tau,t\zeta)
\end{equation}
where
\begin{equation} \label{5.2}
\tilde{p}(t,z,\tilde{\tau},\tilde{\zeta}) \in \ci( \overline{\R}_+  \times \Sigma, \mathfrak{A}^m(B,\g;\R^{1+d}_{\tilde{\tau},\tilde{\zeta}})),\,  \g=(\gamma,\gamma -m, \Theta).       
\end{equation}
Define
\begin{equation} \label{5.3}
M^m_{\mathcal{O}}(B,\g;\R^d_{\zeta})\subset \mathcal{A}(\C,\mathfrak{A}^m(B,\g;\R^d_\zeta)) \ni\,h(v,\zeta)
\end{equation}
such that
\begin{equation}
h(\delta+i\tau,\zeta)\in \mathfrak{A}^m(B,\g;\R^{1+d}_{\tau,\zeta})
\end{equation}
for every $\delta \in\R,$ uniformly in compact $\delta$-intervals.
\begin{Thm}
For $p(t,z,\tau,\zeta)$ as in $\eqref{5.1}$ there exists a $h(t,z,v,\zeta)=\tilde{h}(t,z,v,t\zeta)$ for some $\tilde{h}(t,z,v,\tilde{\zeta})\in\ci(\overline{\R}_+  \times \Sigma,M^m_{\mathcal{O}}(B,\g;\R^d_{\tilde{\zeta}}))$ such that 
\begin{equation} \label{5.4}
\textup{Op}_{t,z}(p)=\textup{op}_M^{\delta-b/2}\textup{Op}_z(h)\,\textup{mod}\, \mathfrak{A}^{-\infty}(\R_+\times\Sigma\times B,\g),
\end{equation}
for $b:= \textup{dim}\,B,$ for all $ \delta \in \R.$
\end{Thm}
Let us now define weighted spaces, first on $B^\wedge= \R_+\times B\, \textup{for compact} \, B\in \mathfrak{M}_1.$ The space $\mathcal{H}^{(s,\gamma,\delta)}(B^\wedge)$ defined to be to be the completion of $\ci(\R_+\times \textup{int}\,B)$ with respect to the norm
\begin{equation} 
{\Big {\{}}(2\pi i)^{-1} \int_{\Gamma_{{\frac{n+1}{2}}-\delta }} \|R^s(\textup{Im}\,v){Mu(z)\|_{\mathcal{W}^{0,\gamma-s}(B)}^2 dv\Big {\}}}^{1/2}
\end{equation}
where $R^s(\lambda)\in\mathfrak{A}^s(B,\g;\R)$ is an order reducing family of edge operators, $b=\textup{dim}\,B,\,\g=(\gamma,\gamma -s, \Theta).$ Moreover, we have the cone spaces $H^{s,\gamma}_{\textup{cone}}(\R\times B)$ obtained as the completion of $\ci(\R \times \textup{int}\,B)$ with respect to the norm $$\Big \{ \int \|\langle t\rangle ^{-s}\textup{Op}_t(p)(\zeta^1)u\|_{\mathcal{W}^{0,\gamma-s}(B)}^2 dt\Big \}^{1/2}$$ for a parameter-dependent elliptic family $p(t,z,\tau,\zeta)= \tilde{p}(z,\langle t\rangle\tau,\langle t\rangle\zeta),$  $\tilde{p}(\tilde{\tau},\tilde{\zeta})\in\mathfrak{A}^s(M,\g;\R^d_{\tilde{\zeta}}),\,\\ \g=(\gamma,\gamma -s, \Theta),$ and $|\zeta^1|$ sufficiently large and fixed. Then we set $H^{s,\gamma}_{\textup{cone}}(B^\wedge):= H^{s,\gamma}_{\textup{cone}}(\R \times B)|_{B^\wedge}.$
Finally for any cut-off function $\omega(t)$ we set 
$$\mathcal{K}^{s,(\gamma,\delta)}(B^\wedge)=\omega \mathcal{H}^{s,(\gamma,\delta)}(B^\wedge) + (1-\omega)H^{s,\gamma}_{\textup{cone}}(B^\wedge),$$
and $\mathcal{K}^{s,(\gamma,\delta);g}(B^\wedge):={\langle t\rangle}^{-g}\mathcal{K}^{s,(\gamma,\delta)}(B^\wedge).$
Similarly as $\eqref{4.3a}$ we form operator-valued amplitude functions
\begin{equation} \label{5.5}
a(z,\zeta) = \omega(t)t^{-m}\op_{M_t}^{\delta - n/2}(h)(z,\zeta)\tilde{\omega}(t) +\textup{smoothing Mellin plus Green symbols}(z,\zeta) 
\end{equation}
belonging to $S^m(\Sigma\times\R^d;\mathcal{K}^{s,(\gamma,\delta)}(B^\wedge),\mathcal{K}^{s-m,(\gamma-m,\delta-m)}(B^\wedge)).$  For $M\in\mathfrak{M}_2$ with the minimal stratum $Z\subset M,$ say, of dimension $d>0,$ we have the space of corner pseudo-differential operators
\begin{equation} 
\mathfrak{A}^m(M,\g)\, \textup{for}\, \g=(\g_1,\g_2),\,\g_1=(\gamma,\gamma -m, \Theta_1),\,\g_2=(\delta,\delta -m, \Theta_2),
\end{equation}
consisting of all $A\in \mathfrak{A}^m(M\setminus Z)$ that are locally near $Z$ of the form $A=\textup{Op}_z(a)\,\textup{mod}\, \mathfrak{A}^{-\infty}(M,\g) $ where $\mathfrak{A}^{-\infty}(M,\g)$ is defined by mapping properties to smooth functions with asymptotics. On $Z$ we have weighted spaces $$\mathcal{W}^{s,(\gamma,\delta)}(M)\subset \mathcal{W}_{\textup{loc}}^{s,\gamma}(M\setminus Z),$$ locally near $Z$ modelled on $\mathcal{W}^s(\R^d,\mathcal{K}^{s,\gamma}(B^\wedge)).$ The principal symbolic structure of an $A\in\mathfrak{A}^m(M,\g)$ is given by $(\sigma(A|_{M\setminus Z}), \sigma_2(A))$ where $\sigma(A|_{M\setminus Z})$ is known from the case $k=1,$ and
$$\sigma_2(A)(z,\zeta)=t^{-m}\op_{M_t}^{\delta - n/2}(h_0)(z,\zeta) $$ for $h_0(t,z,v,\zeta)=\tilde{h}(0,z,v,t\zeta),$ which is a family of operators
$$\sigma_2(A)(z,\zeta): \mathcal{K}^{s,(\gamma,\delta)}(B^\wedge) \rightarrow \mathcal{K}^{s-m,(\gamma-m,\delta-m)}(B^\wedge)$$ for $(z,\zeta)\in T^*Z\setminus0,$ homogeneous in the sense $$\sigma_2(A)(z,\lambda\zeta)=\lambda^m\kappa_{\lambda}\sigma_2(A)(z,\zeta)\kappa_{\lambda}^{-1},\,\lambda>0.$$
If $A$ or $B$ is properly supported (such a property is defined in an analogous manner as in the smooth case) the $AB$ belongs to the corner calculus again, and we have $\sigma(AB)=\sigma(A)\sigma(B)$ with componentwise composition.\\
An operator $A\in\mathfrak{A}^m(M,\g)$ is said to be elliptic if $A|_{M\setminus Z}$ is elliptic in the calculus over $M\setminus Z\in \mathfrak{M}_1,$ and if close to $Z$ the symbolic components $\sigma_0(\cdot,\tilde{\zeta}), \sigma_1(\cdot,\tilde{\zeta})$ with parameter $\tilde{\zeta}\in\R^{\textup{Z}}\setminus \{0\}$ (substituting $t\zeta$) are parameter-dependent elliptic up to $t=0.$ The latter condition concerns an evident generalisation of Definition $\ref{4.6def}$
to the case when we have an extra covariable $\tilde{\zeta}$ which is also involved in edge-degenerate form and where the symbols (apart from a weight factor $t^{-m}$) also depend on $t,$ smoothly up to $t=0.$
\begin{Thm}
Every $A\in\mathfrak{A}^m(M,\g),\,M $ compact, induces continuous operators
\begin{equation} \label{5.6}
 A: \mathcal{W}^{s,(\gamma,\delta)}(M)\rightarrow \mathcal{W}^{s-m,(\gamma-m,\delta-m)}(M)
\end{equation}
 for all $s \in \R$. The operator $\eqref{5.6}$ is compact when $\sigma(A)=0.$
\end{Thm}
\begin{Thm}
Let $A\in\mathfrak{A}^m(M,\g),\,M $ compact; then $A$ is elliptic exactly when $\eqref{5.6}$
 is a Fredholm operator for some $s=s_0 \in \R$. The ellipticity of $A$ entails the Fredholm property of $\eqref{5.6}$ for all $s$. Moreover, if $A$ is elliptic,  $M$ not necessarily compact, the operator $A$ has a parametrix in $\mathfrak{A}^{-m}(M,\g^{-1})$ belonging to $\sigma^{-1}(A)$ \textup{(}with componentwise inverses\textup{)}.
\end{Thm}
\begin{Rem}
The above-mentioned results on parameter-dependent operators in the case $k=1$ have natural analogues for $k=2$. They imply, in particular, the existence of order reducing operators in the calculus.
\end{Rem}
\begin{Exam} \label{5.7ex}
Let $X,Y,$ and $ Z$ be Riemannian manifolds with Riemannian metrics $g_X, g_Y,$ and $g_Z,$ respectively, and form the degenerate metric $$dt^2+t^2(dr^2+r^2g_X+g_Y)+g_Z $$ on the stretched corner $\R_+\times(\R_+\times X\times Y)\times Z.$ Then the associated Laplace-Beltrami operator belongs to the corner calculus for $k=2$ on $(X^{\Delta}\times Y)^{\Delta}\times Z\in\mathfrak{M}_2.$ More generally, considering, for instance, $$M:=(\ldots((X^\Delta\times Y_1)^\Delta\times Y_2)^\Delta\times\ldots Y_{k-1})^\Delta\times Y_k$$ for Riemannian manifolds $X,Y_1,\ldots,Y_k,$ the corner metric
  $$dr_k^2+r_k^2(dr_{k-1}^2+\ldots+(dr_2^2+r_2^2(dr_1^2+r_1^2g_X+g_{Y_1})+g_{Y_2})+\ldots)+g_{Y_k},$$ on $$\R_+\times ( \ldots(\R_+\times (\R_+\times X\times Y_1)\times Y_2)\times\ldots Y_{k-1})\times Y_k$$ gives rise to a Laplace-Beltrami operator belonging to $\textup{Diff}_{\textup{deg}}^2(M).$
\end{Exam}

\section{Concluding remarks} 
 The higher corner calculus that we presented here contains many technicalities that are derived from the program to cover all the substructures sketched in Section 3. A dominating aspect of our theory is to guarantee that the calculus is closed under the construction of parametrices of elliptic elements and that it reflects asymptotics of solutions and elliptic regularity in weighted spaces. What concerns the history of our approach, the above-mentioned information has been integrated from the very beginning, for instance, classical elliptic boundary value problems (BVPs) in the sense of Agmon, Douglis, and Nirenberg \cite{Agmo1}, the theory of pseudo-differential BVPs of Vishik and Eskin \cite{Vivs2}, \cite{Eski2}, the calculus of Boutet de Monvel \cite{Bout1}, details on singular integral operators and operators based on the Mellin transform on the half-axis \cite{Gohb2}, \cite{Eski2}. The symbolic structures have been invented in such a way that vanishing of principal symbols gives rise to compact operators (when the configuration is compact, otherwise after localisation).
Clearly the Fredholm index of elliptic operators has been realised as something invariant under stable homotopies of elliptic principal symbols (through elliptic symbols). At some point the author together with Rempel \cite{Remp7} became aware of the similarity between the boundary symbolic calculus for BVPs without the transmission property, cf. \cite{Eski2}, \cite{Remp1}, and the theory of Kondratiev \cite{Kond1} where $\overline{\R}_+$ is replaced by a cone with a non-trivial base $X.$ The inclusion of edge problems was the next logical step in the development, and, after a preliminary work with Rempel \cite{Remp3}, the paper \cite{Schu32} gave a first systematic edge pseudo-differential calculus. Another step of the development was the paper \cite{Schu29} where the theory has been extended to the case of manifolds with corners (locally modelled on a cone where the base has conical singularities). After that it took some time to develop more technical tools to make the approach really iterative, cf. \cite{Schu25}, \cite{Schu27} (the paper \cite{Schu27} studies singularities modelled on cones where the base has edges).\\
There are many aspects to be deepened and continued in future, for instance, on operator algebras where the symbols do not satisfy an analogue of the Atiyah-Bott condition for the existence of Shapiro-Lopatinskij elliptic edge conditions (for $k\geq 2;$ concerning the case of boundary value problems and edge problems, cf. \cite{Schu37}, and the author's joint papers with Seiler \cite{Schu44}, and \cite{Schu42}), moreover, on the nature of iterated and variable branching asymptotics of solutions, cf. \cite{Schu34}, \cite{Schu36} for the case of boundary value problems, and the joint work with Volpato \cite{Schu60} for the case of edge problems, or the explicit computation of admissible weights or asymptotic data, cf. the author's joint papers with Dines and Liu \cite{Dine3}, \cite{Dine4}, \cite{Liu2} for the case of corner or boundary value problems.\\Different schools on singular analysis apparently emphasise different classes of degenerate operators (in stretched coordinates), and, although there are considerable intersections between the various attempts, it seems that there is no standard terminology on what is a corner-degenerate operator or a corner manifold. Therefore, we point out once again that our calculus is made for corner manifolds that include cones, wedges, cubes, higher polyhedra, etc., embedded in a smooth ambient space, and equipped with the induced (incomplete) corner metrics. Differential operators in the respective stretched coordinates are polynomials in degenerate vector fields of the form 
\begin{align*}
(\partial/\partial x_j)_{j=1,\ldots,n},r_1\partial /\partial r_1,\,&(r_1\partial /\partial y_{1,l})_{l=1,\ldots,q_1},\,r_1r_2\partial /\partial r_2,\,(r_1r_2\partial /\partial y_{2,l})_{l=1,\ldots,q_2},\ldots, \\ &r_1r_2\ldots r_k\partial /\partial r_k,\,(r_1r_2\ldots r_k\partial /\partial y_{k,l})_{l=1,\ldots,q_k},
\end{align*} combined with weight factors $(r_1\ldots r_k)^{-m}$ for operators of order $m,$ and with coefficients that are smooth in all variables up to $r_j=0, j=1,\ldots,k.$ Here $r_j\in \R_+,$ and $(y_{j,l})_{l=1,\ldots,q_j}$ is the variable on a $q_j$-dimensional edge. This is exactly what we obtain as local descriptions of operators $\textup{Diff}_\textup{deg}^m(M),\,M\in\mathfrak{M}_k,$ cf. Section 2, or Example $\ref{5.7ex}$. Moreover, if we are in the situation that an operator is given in a domain with polyhedral boundary, and the respective operator is expressed in Euclidean coordinates in $\R^n$ with smooth coefficients across the boundary (for instance, the standard Laplacian in $\R^n,$ then by repeatedly substituting polar coordinates (according to the order of singularity) we obtain also operators in our class. In such a case it is convenient to formulate everything in the variant of (pseudo-differential) boundary value problems, i.e., to replace the parameter-dependent operators, say, on a closed manifold $X$ (as in Section 4) by the algebra of boundary value problems on $X,$ (now for an X with boundary) with the transmission property
at the smooth faces of the boundary. This aspect is systematically applied in \cite{Kapa10}, \cite{Haru13}, and in numerous other papers mentioned before, jointly with Dines, Liu, Wei, and others. 
%\input{references}

%\addcontentsline{toc}{section}{References}
%\bibliography{/home/schulze/ichbinwichtig/tahoki/bibliography/master}
%\bibliographystyle{/home/schulze/ichbinwichtig/tahoki/bibliography/wileybib}

\addcontentsline{toc}{section}{References}
%\bibliographystyle{wileybib}
%\bibliography{pdeo} 
%\end{document}

% ------------------------------------------------------------------------
\end{document}